\documentclass[a4paper]{amsart} 
\usepackage[ascii]{inputenc} 

\usepackage{amssymb}
\usepackage{mathrsfs}
\usepackage{mathtools}
\usepackage[hidelinks]{hyperref}

\usepackage{xcolor}

\usepackage[shortlabels]{enumitem}
\setlist[enumerate,1]{label={(\Alph*)}}
\setlist[enumerate,2]{label={(\alph*)}}
\setlist[enumerate,3]{label={$\bullet_{\arabic*}$}}

\newenvironment{PROOF}[2][\proofname.]
   {\begin{proof}[#1]\renewcommand{\qedsymbol}{\textsquare$_{\mathrm{#2}}$}}
   {\end{proof}}

\newtheorem{theorem}{Theorem}[section]

\newtheorem{claim}[theorem]{Claim}
\newtheorem{conclusion}[theorem]{Conclusion}

\newtheorem{cruc}[theorem]{Crucial Claim}

\newtheorem{observation}[theorem]{Observation}

\theoremstyle{definition}

\newtheorem{definition}[theorem]{Definition}

\newtheorem{fact}[theorem]{Fact}
\newtheorem{hypothesis}[theorem]{Hypothesis}

\theoremstyle{remark}

\newtheorem{notation}[theorem]{Notation}

\newtheorem{remark}[theorem]{Remark}

\newcommand{\cf}{\mathrm{cf}}

\newcommand{\solv}{\mathrm{solv}}

\newcommand{\dom}{\mathrm{dom}}

\newcommand{\GCH}{\mathrm{GCH}}

\newcommand{\Ord}{\mathrm{Ord}}
\newcommand{\otp}{\mathrm{otp}}

\newcommand{\Wilog}{\ensuremath{\text{Without loss of generality}}}

\newcommand{\bfA}{\mathbf{A}}
\newcommand{\bfa}{\mathbf{a}}

\newcommand{\bfb}{\mathbf{b}}
\newcommand{\bfc}{\mathbf{c}}

\newcommand{\bfG}{\mathbf{G}}

\newcommand{\bfn}{\mathbf{n}}
\newcommand{\bfp}{\mathbf{p}}

\newcommand{\bfQ}{\mathbf{Q}}
\newcommand{\bfq}{\mathbf{q}}

\newcommand{\bfs}{\mathbf{s}}

\newcommand{\bfV}{\mathbf{V}}

\newcommand{\bbL}{\mathbb{L}}

\newcommand{\bbP}{\mathbb{P}}
\newcommand{\bbQ}{\mathbb{Q}}
\newcommand{\bbR}{\mathbb{R}}

\newcommand{\cB}{\mathscr{B}}

\newcommand{\cH}{\mathscr{H}}
\newcommand{\cI}{\mathscr{I}}

\newcommand{\cS}{\mathscr{S}}

\newcommand{\cU}{\mathscr{U}}

\newcommand{\gC}{\mathfrak{C}}

\newcommand{\varp}{\varepsilon}

\newcommand{\rest}{\restriction}

\newcount\skewfactor
\def\mathunderaccent#1#2 {\let\theaccent#1\skewfactor#2
\mathpalette\putaccentunder}
\def\putaccentunder#1#2{\oalign{$#1#2$\crcr\hidewidth
\vbox to.2ex{\hbox{$#1\skew\skewfactor\theaccent{}$}\vss}\hidewidth}}
\def\name{\mathunderaccent\tilde-3 }

\newbox\noforkbox \newdimen\forklinewidth
\setbox0\hbox{$\textstyle\bigcup$}
\setbox1\hbox to \wd0{\hfil\vrule width \forklinewidth depth \dp0
                        height \ht0 \hfil}
\setbox\noforkbox\hbox{\box1\box0\relax}
\def\unionstick{\mathop{\copy\noforkbox}\limits}
\def\nonfork#1#2_#3{#1\unionstick_{\textstyle #3}#2}
\def\nonforkin#1#2_#3^#4{#1\unionstick_{\textstyle #3}^{\textstyle
    #4}#2}
\setbox0\hbox{$\textstyle\bigcup$}
\setbox1\hbox to \wd0{\hfil{\sl /\/}\hfil}
\setbox2\hbox to \wd0{\hfil\vrule height \ht0 depth \dp0 width
                                \forklinewidth\hfil}
\newbox\doesforkbox
\setbox\doesforkbox\hbox{\box1\box0\relax}
\def\nunionstick{\mathop{\copy\doesforkbox}\limits}

\def\fork#1#2_#3{#1\nunionstick_{\textstyle #3}#2}
\def\forkin#1#2_#3^#4{#1\nunionstick_{\textstyle #3}^{\textstyle
    #4}#2}

\newcommand{\stickT}{%
\setbox255=\hbox{\raise1ex\hbox{$\hspace{0.2pt}\,\bullet\,$}}
\mathord{\rlap{\hbox to\wd255{\hss\hbox{$|$}\hss}}
\box255}
}
\newcommand{\stickS}{%
\setbox255=\hbox{\raise0.6ex\hbox{$\scriptstyle\bullet$}}
\mathord{\rlap{\hbox to\wd255{\hss\hbox{$\scriptstyle|$}\hss}}
\box255}
}

\usepackage[hidelinks]{hyperref}

\author[S. Shelah]{Saharon Shelah}
\address{Einstein Institute of Mathematics,
The Hebrew University of Jerusalem,
9190401, Jerusalem, Israel; and\\
Department of Mathematics,
Rutgers University,
Piscataway, NJ 08854-8019, USA}
\urladdr{https://shelah.logic.at/}

\thanks{First typed: January 7, 2024.
The author would like to thank the Israel Science Foundation for partial support of this research by grant 2320/23. The author would like to thank Craig Falls for generously funding typing services, and the typist for his careful and beautiful work. 
Submitted to 
Yair Hayut in February 2025 
for the IJM volume dedicated to Menachem 
Magidor.  References like [Sh:950, Th0.2=Ly5] mean that the internal label of Th0.2 is y5 in Sh:950.
The reader should note that the version in my website is usually more up-to-date than the one in arXiv. This is publication  number 1258 on Saharon Shelah's list.}

\makeatletter
\@namedef{subjclassname@2020}{\textup{2020} Mathematics Subject Classification}
\makeatother

\date{December 24, 2025} 

\title{Consistency of square bracket partition relation}

\begin{document}
\makeatletter\def\shfiuwefootnote{\gdef\@thefnmark{}\@footnotetext}\makeatother\shfiuwefootnote{Version 2026-01-06. See \url{https://shelah.logic.at/papers/1258/} for possible updates.}
\begin{abstract}
    Characteristic earlier results were of the form ${\mathrm{CON}}(2^{\aleph_0} \to [\lambda]_{n,2}^2)$, with $2^{\aleph_0}$ 
    an ex-large  
    cardinal,  in  
    the best case 
    the first weakly Mahlo cardinal.

    Characteristic new results are ${\mathrm{CON}} \big( (2^{\aleph_0} = \aleph_m) +
    \aleph_{\ell} \to [\aleph_k]_{n,2}^2 \big)$, 
     for suitable 
    $k < \ell < m$. So we improve in three respects: the continuum may be small (e.g. not a weakly Mahlo), we use no large cardinal, and the cardinals 
    $ \lambda $  
    involved are $<2^{\aleph_0}$ after the forcing. 

\end{abstract}
\maketitle

\setcounter{section}{-1}

\section{Introduction}

In their seminal list of problems \cite{EH}, Erd\"os and Hajnal posed the question ($15$(a)): does $2^{\aleph_{0}} \not\to [\aleph_{1}]_{3}^{2}$?  
Recently, Komj\'ath~\cite{Kom25} provided a comprehensive update on this topic.

We continue here works which start with the problem above:\cite[\S 2]{Sh:276}, \cite{Sh:288}, \cite{Sh:289}, \cite{Sh:473}~\cite{Sh:481}, \cite{Sh:546} and the work with Rabus~\cite{Sh:585}, but we try to be self-contained. 

The simplest case of our result is (recall \ref{x8} below): 

\begin{theorem}\label{x2}
    Assume $\GCH$ for transparency. \underline{Then} for some ccc forcing notion of cardinality $\aleph_{6}$
    in the universe $\bfV^{\bbP}$, we have $2^{\aleph_{0}} = \aleph_{6}$ 
      and for any $n \geq 3$, $\aleph_{5}
     \to [\aleph_{2}]_{n, 2}^{2}$. 
\end{theorem}

\begin{PROOF}{\ref{x2}}  
Choose  $ (\mu, \theta, \partial , \lambda )$
as 
$ ( \aleph _6,
 \aleph _5, 
 {\aleph_2}, 
 {\aleph_0} )$
and apply Theorem  \ref{x5}  and Fact \ref{d26} with $\partial_{0} = \aleph_{1}$. 
\end{PROOF}

For Hypothesis \ref{d2}, the main case is:

\begin{theorem}\label{x5}
    Assume $\lambda = \lambda^{< \lambda} < \partial < \theta < \mu = \mu^{\theta}$,
$\partial = \partial^{\lambda}$ and $2^{\partial^{+ \ell}} = 
   \partial^{+ \ell + 1}$ for $\ell = 0, 1, 2$ and 
    $\partial^{+4} \leq \theta$.  
     \underline{Then} for some $\lambda^{+}$-cc, $(< \lambda)$-complete forcing notion $\bbP$ 
    of cardinality $ \mu $  
    (so the forcing does not collapse any cardinal
    and preserves 
    cardinal arithmetic outside $ [\lambda, \mu)$),  
    in the universe $\bfV^{\bbP}$ we have, 
    $ 2^ \lambda = \mu $  and 
    for every $\sigma < \lambda$, $\theta \to [\partial]_{\sigma, 2}^{2}$  
\end{theorem}

\begin{PROOF}{\ref{d23}}
  All this paper is dedicated to proving
  this theorem.
  Pedantically, choose $\partial = \kappa^{+}$, notice that Hypothesis \ref{d2}  holds 
  (by Fact \ref{d26}) so we can apply
  Conclusion \ref{d23}.
\end{PROOF}

We may weaken $\mu = \mu^{\theta}$ to $\mu = \mu^{\partial}$ and replace $\partial = \kappa$ by $\partial$ being a suitable limit cardinal. 


Recall, 

\begin{definition}\label{x8}
    For   
    possibly finite  
    cardinals $\theta, \partial, \sigma$ and $\kappa$, let $\theta \to [\partial]_{\sigma, \kappa}^2 
    $ mean: 

    \begin{itemize} 
        \item if  $ \bfc $ 
        is a function from $[\theta]^{2} \coloneqq \{ u \subseteq \theta \colon \vert u \vert = 2 
        \}$ into $\sigma$, \underline{then} there exists some subset $\cU$ of $\theta$ of cardinality $\partial$ such that $\{ \bfc(u) \colon u \in [\cU]^{2}  \}$ has at most $\kappa$-many members. 
    \end{itemize}
\end{definition}

We thank Yair Hayut and the referee for many helpful comments. 

\subsection{Preliminaries}\label{0B}

\begin{notation}\label{z2}
    \ 
    
    1) $\mathrm{cof}(\delta)$ is the class of ordinals of cofinality $\cf(\delta)$.

    2) For a set $x$, let $\mathrm{trcl}(x)$ be the transitive closure of $x$, that is, the minimal set $Y$ such that $x \in Y$ and $(\forall y)(y \in Y \Rightarrow y \subseteq Y)$. 

    3) Let $\cH(\lambda) \coloneqq \{ x \colon \vert \mathrm{trcl}(x) \vert < \lambda \}$. 

    4) Let $\mathrm{trcl}_{\Ord}(x)$ be defined similarly to $\mathrm{trcl}(x)$ considering ordinals as atoms (= elements), equivalently, the minimal set $Y$ such that $x \in Y$ and 
    \[
    (\forall y)[y \in Y \wedge ( \text{if } y \text{ is not an ordinal, then } y \subseteq Y)].
    \]

    5) Let $\cH_{< \kappa}(x) = \{ x \colon \mathrm{trcl}_{\Ord}(x) \subseteq \cH(x) \text{ but has cardinality} < \kappa \}$.
\end{notation}

\begin{notation}\label{z5}\ 
    \
    
    (1) $\bbP$, $\bbQ$ and $\bbR$ are forcing notions. 

    (2) $p, q, r$ called \emph{conditions} are members of a forcing notion. 

    (3) $\bfq$ is as in Definition~\ref{d8}, some kind of $(< \lambda)$-support iterated forcing with extra information. 
\end{notation}

\begin{notation}\label{z8}
    We may write e.g. $N[\bfq, \beta, u]$ instead $N_{\bfq, \beta, u}$ to help with sub-scripts (or super-script). 
\end{notation}

\begin{definition}\label{z11}
    Let $\theta, \partial
    , \kappa  
    $ and $\lambda$ be infinite cardinals. We say that $\theta \to_{\mathrm{sq}} (\partial)_{\kappa 
    }^{\lambda 
    , 2}$ \underline{when} $\theta > \partial \geq 
    \kappa \geq \lambda $ 
    and: 

    \begin{enumerate}
        \item[$\boxplus$] If (a) then  (b), where: 

        \begin{enumerate}
            \item[(a)] $\cB$ is an expansion of $(\cH(\chi), \in, <_{\ast})$, where $<_{\ast}$ is a well-ordering of $\cH(\chi)$, $\chi > \theta$, and its vocabulary $\tau_{\cB}$ has cardinality $\leq \lambda$. 

            \item[(b)] There is a tuple  $\bfs = (\cU, \bar{N}, \bar{\pi})$ solving $\bfp = (\mu, \theta, \partial, \kappa, \lambda, \cB)$, 
               which means: 

            \begin{enumerate}
                \item[$\boxplus_{\bfp, \bfs}$] for $u, v \in [
                {\mathscr U } 
                ]^{\leq 2}$,
        
                \begin{itemize}
                    \item[$\bullet_{1}$] $\bar{N} = \langle N_{u} \colon u \in [\cU]^{\leq 2} \rangle$, 
                    
                    \item[$\bullet_{2}$] $\cU \subseteq \theta$ is such that $\otp(\cU) = \partial$, 
        
                    \item[$\bullet_{3}$] $N_{u} \prec \cB$, $[N_{u}]^{< 
                    \lambda 
                    } \subseteq N_{u}$,  
        
                    \item[$\bullet_{4}$] $\varp[\bfs] \coloneqq \min(\cU
                    )$, 
        
                    \item[$\bullet_{5}$] $N_{u} \cap \cU = u$,
        
                    \item[$\bullet_{6}$] $\Vert N_{u} \Vert = 
                     \kappa 
                    $ and 
                    $  \kappa 
                    + 1 \subseteq N_{u}$, 
        
                    \item[$\bullet_{7}$] $N_{u} \cap N_{v} \prec N_{u \cap v}$,
        
                    \item[$\bullet_{8}$]  $\bar{\pi} = \langle \pi_{u, v} \colon u, v \in [\cU]^{\leq 2 
                    }$ and $\vert u \vert = \vert v \vert \rangle$ 
                    such that if $\vert u \vert = \vert v \vert 
                           $, 
                        then $\pi_{u, v}$ is an isomorphism from $N_{v}$ onto $N_{u}$ mapping $v$ onto $u$, 
        
                    \item[$\bullet_{9}$] if $u_{1} \subseteq u_{2}$ and $v_{1} \subseteq v_{2}$ all from $[\cU]^{\leq 2}$ and $\vert u_{2} \vert = \vert v_{2} \vert$,  
                    $\pi_{u_{2}, v_{2}}
                    ''(v_{1}) = u_{1}$ \underline{then} 
                    $\pi_{u_{1}, v_{1}}$,
                     $\pi_{u_{2}, v_{2}}$ are
                      compatible functions\footnote{So e.g. it follows that: if $\zeta_{1}, \zeta_{2} \in \cU$ then $\pi_{\{ \zeta_{1} \}, \{ \zeta_{2} \}} \rest (N_{\emptyset} \cap N_{\{ \zeta_{2} \}})$ is the identity map.
                      }, 

            \item[$ \bullet _{10}$] 
            for $ {\ell} = 1,2 $, the sets
            $ N_u \cap \partial $  for $ u \in [{\mathscr U } ]^{{\ell} }$                
            are pairwise equal
        \footnote{Note that $ \partial $ has two distinct roles: the size of $ {\mathscr U } $ and the restriction 
 on $ N_ u \cap \partial $. We may separate.} 
        and included in $ N_ \emptyset $.
                \end{itemize}
            \end{enumerate}
        \end{enumerate}
    \end{enumerate}
\end{definition}

\begin{observation}\label{z14}
If $ \bar{N} = \langle N_u: u \in [ {\mathscr U } ]^{\le 2} \rangle  $ 
satisfies \ref{z11}(b)$\bullet_{1} + \bullet_{7}$,  then:
\begin{enumerate} 
\item[$(\ast)$]  For every   $ x \in \cup   \{  N_ u:
  u \in [ {\mathscr U } ]^{\le 2} \} $ 
  the set $ \{   u \in [{\mathscr U } ]^{\le 2} : x \in N_ u \} $ 
  has one of the following forms: 
  \begin{enumerate} 
  \item[(a)] $ \{ u \} $  for some $ u \in [ {\mathscr U } ]^{2}$,
  
  \item[(b)] $ \{ \zeta \} $  for some $ \zeta \in \mathscr{U}$,
  
  \item[(c)] $ \{ \{ \zeta \} \}  \cup \{  \{ \varepsilon , \zeta \}: \varepsilon \in {\mathscr U } \cap  \zeta \}$ for some $ \zeta \in  \mathscr{U},$
  
   \item[(d)]   $  \{ \{ \zeta \} \}
     \cup \{  \{ \zeta, \xi  \}: \xi \in {\mathscr U } \setminus  (\zeta +1) \} $  for some $ \zeta \in \mathscr{U}$,
     
    \item[(e)]  $ \{  \emptyset \}$,
    
    \item[(f)] $ \{  \emptyset \} \cup  \{\{ \zeta \} \colon \zeta \in \mathscr{U}\}$,
    
    \item[(g)] $ \{  \emptyset \} \cup  \{\{ \zeta \} \colon \zeta \in {\mathscr U }    \}  \cup  \{ \{\varepsilon, \zeta  \} \colon  \varepsilon < \zeta  \text{ are from } {\mathscr U } \} $.
    \end{enumerate} 
    \end{enumerate} 
\end{observation} 

\section{The forcing}\label{1}

Our aim here is to prove the consistency of the following configuration: 
\[
2 < \sigma < \lambda = \lambda^{< \lambda} < \partial = \partial^{< \lambda} < \theta < \mu = \mu^{\theta} = 2^{\lambda},
\]
and having $\theta \to [\partial]_{\sigma, 2}^2$. 


A continuation is in preparation \cite{Sh:F2407}, aiming to further develop the directions explored here, particularly for the case of superscript 
$\bfn > 2$, as dealt within~\cite{Sh:288}. We also show there that we can weaken the
requirements on the cardinals and have more pairs.  


\begin{hypothesis}\label{d2}
    The parameter $\bfp = ( \mu,  
    \theta, \partial, 
    \lambda, \lambda 
         , \cB )$   
      consists of the following: 

    \begin{enumerate}
        \item[(a)] $\lambda = \lambda^{< \lambda} < \partial < \theta < \mu = \mu^{\theta}$, 

        \item[(b)] $\theta \to_{\mathrm{sq}}
        (\partial)_{
        \lambda } 
        ^{ \lambda, 2}$ (see Definition~\ref{z11}, a variant of \cite[2.1]{Sh:289}); 
        in our case using 
        $ \lambda $  
        twice in intentional.  

        \item[(c)] $\sigma$ will vary on the cardinal numbers from 
        $ (2, \lambda)$ and the 
        ``nice" 
        $\mu$-s are such that $\gamma < \mu \Rightarrow \vert \gamma \vert^{\theta} < \mu$.  

        \item[(d)] 

        \begin{itemize}
            \item $\chi$ is e.g. $\beth_{2}(\mu)^{+}$,
            

            \item let $\cB$ be an expansion of  $(\cH(\chi), \in, <_{\chi}^{\ast})$ with vocabulary of cardinality $\lambda$ such that for any finite set $u \subseteq \cH(\chi)$, the Skolem hull of $u$ $N_{u}  \coloneqq \mathrm{Sk}(u, \gC_{\ast})$ is of cardinality $\lambda$ and $\vert N_{u} \vert^{< \lambda} \subseteq N$.
        \end{itemize}
    \end{enumerate}
\end{hypothesis}

We intend to use $({<} \, \lambda)$-support iterated forcing of quite a special kind but first, we define the iterand. 

\begin{definition}\label{d5}\ 

     (1) Let $\bfA$ be the set of objects $\bfa$ consisting of: 

    \begin{enumerate}
        \item[(a)] 

        \begin{itemize}
            \item $\gamma  < 
     \mu$ and $\sigma \in (2, \lambda)$, 

            \item $\bbP$ is a forcing notion such that: 
            $$
            p \in \bbP \Rightarrow \dom(p) \in [\gamma]^{< \lambda} \wedge (\forall \alpha \in \dom(p))(p(\alpha) \in [\lambda \cup \gamma]^{< \lambda}),
            $$ 

            \item $\bbP$ is $\lambda^{+}$-cc and $(< \lambda)$-complete, 

            \item the order $\leq_{\bbP}$ is: $p \leq_{\bbP} q$ iff: 
            $$
            \dom(p) \subseteq \dom(q) \wedge (\forall \alpha \in \dom(p))[p(\alpha) \subseteq q(\alpha)],
            $$ 
        \end{itemize}

        \item[(b)] 

        \begin{itemize} 
            \item $\name{\bfc}$ is a $\bbP$-name of a function from $[\theta]^{2}$ to $\sigma$, (we may write $\name{\bfc}(\alpha, \beta)$ instead $\name{\bfc}
            ( 
            \{ \alpha, \beta \}
            )$  
            for $\alpha \neq \beta < \theta$).
        \end{itemize}

        \item[(c)]  We have
        $(\cU, \bar{N}, \bar{  \pi} )$ solving $\bfp = (\mu, \theta, \partial, \lambda , 
        \lambda, \cB)$,  (with 
        $ \cB $ as in  
        Definition \ref{z11}$\boxplus$(b)
          and 
        \ref{d2})  
         such that $\bbP, \name{\bfc} \in N_{u}$ for every $u \in [\cU]^{\leq 2 }$.
    \end{enumerate}

    (1A) In the context of Definition~\ref{d5}(1), $\bfa = (\gamma, \bbP, \name{\bfc}, \cU, \bar{N}, \bar{\pi}) = (\gamma_{\bfa}, ...)$, so e.g. $N_{\bfa, u} = N_{u}$.

    (2) We say that the 
  pair
    $(p, \bar{\iota}
    )$ is a \emph{solution} of $\bfa \in \bfA$, and write $(\bfa, p, \bar{\iota}
    ) \in \bfA^{+}$, \underline{when}:  

    \begin{enumerate}
        \item[(a)] $\bar{\iota} = (\iota_{1}, \iota_{2}) \in \sigma \times \sigma$,

        \item[(b)] $p \in \bbP_{\bfa} \cap N_{\bfa, \{ \varp[\bfa] \}}$,
    recalling $ \varepsilon ( \mathbf{a} ) = \min ({\mathscr U } ) $,  

        \item[(c)] if $p \leq q \in \bbP_{\bfa} \cap N_{\bfa, \{ \varp[\bfa] \}}$ and $\zeta_{1} < \zeta_{2}$ are from $\cU$ \underline{then} there are $q_{1}, q_{2}, r_{1}, r_{2}$ such that for $\ell = 1, 2$, we have: 

        \begin{itemize}
            \item[$\bullet_{0}$] $q \leq_{\bbP_{\bfa}} q_{\ell}$,
            
            \item[$\bullet_{1}$] $q_{\ell} \in \bbP_{\bfa} \cap N_{\bfa, \{ \varp[\bfa] \}}$ and $q_{1} \rest (N_{\bfa, \emptyset} \cap \gamma_{\bfa}) = q_{2} \rest (N_{\bfa, \emptyset} \cap \gamma_{\bfa})$,

            \item[$\bullet_{2}$] $r_{\ell} \in \bbP_{\bfa} \cap N_{\bfa, \{ \zeta_{1}, \zeta_{2} \}}$,

            \item[$\bullet_{3}$] $r_{\ell} \Vdash$``$\name{\bfc}(\zeta_{1}, \zeta_{2}) = \iota_{\bfa, \ell}$'',

            \item[$\bullet_{4}$] $r_{\ell} \rest N_{\bfa, \{ \zeta_{1} \}}$ is $\leq_{\bbP_{\bfa}}$-below $\pi_{\{ \zeta_{1} \}, \{ \varp[\bfa] \}}^{\bfa}(q_{\ell})$, 

            \item[$\bullet_{5}$] $r_{\ell} \rest N_{\bfa, \{ \zeta_{2} \}}$ is $\leq_{\bbP_{\bfa}}$-below $\pi_{\{ \zeta_{2} \}, \{ \varp[\bfa] \}}^{\bfa}(q   
            _{3 - \ell})$. 
        \end{itemize}



    \end{enumerate}

(3) If $\bfb = (\bfa, p, \bar{\iota}) \in \bfA^{+}$ 
then let $\name{\bbQ}_{\bfb}$ be the $\bbP$-name of the following forcing notion: 

    \begin{enumerate}
        \item[$(\ast)$] For $\bfG \subseteq \bbP$ generic over $\bfV$, 

        \begin{enumerate}
            \item[(a)] the set of elements of $\bbQ_{\bfb} = \name{\bbQ}_{\bfb}[\bfG]$ is:
                \begin{multline*}
                 \Bigl\{ u \in [\cU]^{< \lambda} \colon \text{ if } \zeta_{1} < \zeta_{2} \text{ in } \cU, \text{ then } \name{\bfc}
                \{ \zeta_{1}, \zeta_{2} \}[\bfG] \in \{ \iota_{1}, \iota_{2} \}, \text{ moreover}\\
                \text{for some } q_{1}, q_{2}, r_{1}, r_{2} \text{ as in Definition~\ref{d5}(1)(c)$(\bullet_{1})$-$(\bullet_{5})$}, \text{ we have }  r_{1} \in \bfG \text{ or } r_{2} \in \bfG \Bigr \},
                \end{multline*}

            \item[(b)] the order of $\name{\bbQ}_{\bfb}[\bfG]$ is 
inclusion, 

            \item[(c)] the generic is $\name{\mathcal{V}}_{\bfb} = \bigcup \name{\bfG}_{\name{\bbQ}_{\bfb}}$. 
        \end{enumerate}
    \end{enumerate}
\end{definition}

\begin{definition}\label{d8}\

    (1) Let $\bfQ \coloneqq \bfQ_{\bfp}$ be the class of $\bfq$ which consist of (below, $\alpha \leq \lg(\bfq)$ and $\beta < \lg(\bfq)$ and e.g. $\bbP_{\alpha} = \bbP_{\bfq, \alpha}$): 

    \begin{enumerate}
        \item[(a)] $\lg(\bfq)$ is an ordinal $\leq \mu$, 

        \item[(b)] $\langle \bbP_{\alpha}, \name{\bbQ}_{\beta} \colon \alpha \leq \lg(\bfq), \beta < \lg(\bfq) \rangle$ is a $({<} \, \lambda)$-support iteration, 

        \item[(c)] $\bbP_{\beta}$ satisfies the $\lambda^{+}$-cc, 

        \item[(d)] $\name{\bbQ}_{\beta}$ is $\name{\bbQ}_{\bfb_{\beta}}$, where: 

        \begin{itemize}
            \item[$\bullet_{1}$] $\bfb_{\beta} \coloneqq (\bfa_{\beta}, p_{\beta}^{\ast}, \bar{\iota}_{\beta}^{\ast}
            ) \in \bfA^{+}$, 

            \item[$\bullet_{2}$] $\bfa_{\beta} \coloneqq (\gamma_{\beta}, \bbP_{\beta}^{\bullet}, \name{\bfc}_{\beta}, \cU_{\beta}, \bar{N}_{\beta}, \bar{\pi}_{\beta}) \in \bfA$,

            \item[$\bullet_{3}$] $\bbP_{\beta}^{\bullet}$ is equal to $\bbP_{\xi(\beta)}'$ for some $\xi(\beta) = \xi_{\bfq}(\beta) \leq \beta$
            (on $\bbP_{\beta}'$, see below),

          \item[$ \bullet _4$] 
          The sequence 
          $ \langle (\mathbb{P} _ 
             \gamma, \mathbb{P} '_\gamma , \mathbf{a} _\gamma , \mathbf{b} _\gamma , 
          \xi 
          ( \gamma )  ) :
          \gamma < \beta \rangle $
          belongs to $ N_{\beta , u}$ 
          for every $ u \in [{\mathscr U }_ \beta  
            ]^{\le 2}$. 

          \item[$ \bullet _5$] 
Let $ {\mathscr W } _ \beta = \bigcup \{ N_{\beta, u} \cap  \beta \colon 
u \in [{\mathscr U }
_ \beta  
]^{\le 2 } \}, $   

\item[$\bullet_{6}$]
we 
   \footnote{
      Why?  
      By \ref{z11}(b)$\bullet_{10}.$   
   }
    have: for every $ \gamma \in \mathscr{W}_{\beta}$ the set $\mathscr{W}_{\beta} \cap \mathscr{W}_{\gamma}$ has cardinality $\leq \lambda$, 

    \item[$\bullet_{7}$] For every $\gamma \in \mathscr{W}_{\beta}$, there is  $ u  
    = u_{\beta , \gamma }  
    \in [ {\mathscr U } _ \beta  
    ]^{\le 2
      } $ 
    such that ${\mathscr W } _\beta \cap {\mathscr W } _ \gamma  \subseteq 
    N_{\beta, u}$ and without loss of generality $u$ is minimal with this property.  
        \end{itemize}

        \item[(e)] $\bbP_{\alpha}'$ is a dense subset of $\bbP_{\alpha}$, where, 

        \begin{itemize}
            \item $\bbP_{\alpha}'$ is $\bbP_{\alpha}$ restricted to the set of conditions $p \in \bbP_{\alpha}$ such that:
            
            if $\beta \in \dom(p)$ then $p(\beta)$ is a member of $\bfV$ (not just a $\bbP_{\alpha}$-name) and if $\zeta_{1} < \zeta_{2}$ 
        are in $ p(\beta) \subseteq 
        \cU_{\beta}$, then there are $q_{1}, q_{2}$, $r_{1}$, $r_{2}$ as in Definition~\ref{d5}(2)(c)$(\bullet_{1})$-$(\bullet_{5})$ with  $\bfa_{\beta}, \bfb_{\beta}$ here standing for $\bfa, \bfb  
            $ there and 
            $$
            \bigvee_{\ell = 1}^{2}(\forall \gamma \in \dom(r_{\ell}))[ \gamma \in \dom(p) \wedge r_{\ell}(\gamma) \subseteq p(\gamma)].
            $$ 
        \end{itemize}

        \item[(f)] $\gamma_{\bfq} \coloneqq \gamma(\bfq) \coloneqq \sup \{ \gamma_{\bfq, \beta} \colon \beta < \lg(\bfq) \}$, so $\bbP_{\gamma(\bfq)}' \subseteq \cH_{< \lambda}(\gamma_{\bfq})$; let $\bbP_{\bfq} \coloneqq \bbP_{\lg(\bfq)}$ and $\bbP' 
        _ \bfq 
        \coloneqq \bbP_{\lg(\bfq)}'$. 


    \end{enumerate}

    (1A) We may write either $\bbP_{\bfq, \alpha}$  or $\bbP_{\alpha}$ whenever $\bfq$ is clear and $(
    {\iota}_{\bfq, \beta, 1}, 
    {\iota}_{\bfq, \beta, 2})$ is $\bar{\iota}_{\bfb  _ 
    \beta}$. 

    (2) Let $\leq_{\bfp}$ be the following two-place relation on $\bfQ_{\bfp}$: 
    $$
    \bfq_{1} \leq_{\bfp} \bfq_{2} \text{ \underline{iff} } \bfq_{1} = \bfq_{2} \rest \lg(\bfq_{1}), \text{ see below. }
    $$

    (3) For $\bfq_{2} \in \bfQ_{\bfp}$ and $\alpha_{\ast} \leq \lg(\bfq_{2})$, we define $\bfq_{1} \coloneqq \bfq_{2} \rest \alpha_{\ast}$ by: 

    \begin{enumerate}
        \item[(a)] $\lg(\bfq_{1}) = \alpha_{\ast}$, 

        \item[(b)] $( \bbP_{\bfq_{1}, \alpha}, \bbP_{\bfq_{1}, \alpha}' ) = (\bbP_{\bfq_{2}, \alpha}, \bbP_{\bfq_{2}, \alpha}')$ for $\alpha \leq \alpha_{\ast}$,

        \item[(c)] $(\name{\bbQ}_{\bfq_{1}, \beta}, \bfb_{\bfq_{1}, \beta}
     , \xi  _{\mathbf{q} _1}
            ( \beta )  
        ) = (\name{\bbQ}_{\bfq_{2}, \beta}, \bfb_{\bfq_{2}, \beta}
    , \xi _{\mathbf{q} _2} (\beta     
        ))$ for $\beta < \alpha_{\ast}$. 
    \end{enumerate}

    (4) We say that two conditions $p, q \in \bbP_{\alpha}'$ are \emph{isomorphic}, \underline{when}: 

    \begin{enumerate}
        \item[(a)] $\otp(\dom(p)) = \otp(\dom(q))$, and

        \item[(b)] if $\beta \in \dom(p) \cap \dom(q)$ then: 

        \begin{itemize}
            \item[$\bullet_{1}$] $\otp(p(\beta)) = \otp(q(\beta))$,

            \item[$\bullet_{2}$] if $\varp \in  
            p(\beta) \cap q(\beta)$ then $\otp(\varp \cap p(\beta)) = \otp(\varp \cap q (\beta))$, 

            \item[$\bullet_{3}$] if $\varp \in p(\beta), \zeta \in q(\beta)$ and $\otp(\varp \cap p(\beta)) = \otp(\zeta \cap q(\beta))$ then: 
            $$
            \pi_{\beta, \{ \zeta \}, \{ \varp \}} (p \rest N_{\beta, \{ \varp \}}) =
    q \rest N_{\beta, \{ \zeta \}}.
            $$
    
        \item[$ \bullet _4$] if $ \varepsilon < 
      \varepsilon _1 $ belong 
      to   
        $ p(\beta ) , 
        \zeta < \zeta_1 $ 
          belong 
          to 
          $ q(\beta )$, 
         $\otp(\varp \cap p(\beta)) = \otp(\zeta \cap q(\beta))$   and 
          $\otp(\varp _1 \cap p(\beta)) = \otp(\zeta _1 \cap q(\beta))$ 
      then:
      $$ \pi _{\beta , \{  \zeta, \zeta_1\} , \{  \varepsilon, \varepsilon_1\}  }  
      (p \upharpoonright N_{\beta , \{ \varepsilon, \varepsilon _1\} } )=
      q \upharpoonright N_{\beta, \{ \zeta , \zeta_1 \}  
      }.
      $$
        \end{itemize}
    \end{enumerate}
\end{definition}


\begin{remark}\label{d9}
    If we prefer in clause (d)$~(\bullet_{3})$ of Definition \ref{d8}~(1) to have $\xi(\beta) = \beta$, i.e., $\bbP_{\beta}^{\bullet} = \bbP_{\beta}'$, we need to add, e.g. 
    ``$\mu$ is regular and  
    e.g.  
    use a preliminary forcing
    $( \{ \bfq \in \bfQ_{\bfp} \colon \lg(\bfq) < \mu \}, \lhd )$''. 
\end{remark}

\begin{claim}\label{d11}
    \ 

    (0) For $\bfq \in \bfQ_{\bfp}$, we have: $\bbP_{\bfq}' \models$``$p \leq q$'' iff $\{p, q \} \subseteq \bbP_{\bfq}'$, $\dom(p) \subseteq \dom(q)$, and $\beta \in \dom(p) \Rightarrow p(\beta) \subseteq q(\beta)$. 
    
    (1) For $\bfq \in \bfQ_{\bfp}$, any increasing sequence of members of length $< \lambda$ of $\bbP_{\bfq}'$ has a lub, in fact, if $\delta < \lambda$, $\bar{p} = \langle p_{i} \colon i < \delta \rangle \in {}^{\delta} (\bbP_{\bfq}')$ is increasing, \underline{then} the following $p \in \bbP_{\bfq}'$  is a lub of $\bar{p}$; defined by: $\dom(p) = \bigcup \{ \dom(p_{i}) \colon i < \delta \}$, and if $\beta \in \dom(p)$ then 
    $$
    p(\beta) = \bigcup \left \{ p_{i}(\beta) \colon i < \delta \text{ and } \beta \in \dom(p_{i}) \right\}.
    $$ We denote this $p$ by $\lim(\bar{p})$. 

    (2) For $\bfq \in \bfQ_{\bfp}$, we have: 

    \begin{itemize}
        \item[$\bullet$] $p \in \bbP_{\bfq}'$ \underline{iff}: 

        \begin{enumerate}
            \item[(a)] $p$ is a function with domain $\in [\lg(\bfq)]^{< \lambda}$, 

            \item[(b)] if $\beta \in \dom(p)$ then $p(\beta)$ belongs to $[\cU_{\beta}]^{< \lambda}$.

            \item[(c)] If $\beta \in \dom(p)$ and $(\iota_{1}, \iota_{2}) = (\iota_{\bfq, \beta, 1}, \iota_{\bfq, \beta, 2})$ then  for every $\zeta_{1} < \zeta_{2}$ from $p(\beta)$, $(p \rest \beta) \, {\rest} \, N_{\bfq, \beta, \{ \zeta_{1}, \zeta_{2} \}} \Vdash_{\bbP_{\bfq, \beta}}$``$\name{\bfc} \{ \zeta_{1}, \zeta_{2} \} \in \{ \iota_{1}, \iota_{2} \}$''. Moreover, there are $q_{1}$, $q_{2}$, $r_{1}$, $r_{2}$ as in Definition~\ref{d5}(2)(c)$(\bullet_{1})$-$(\bullet_{5})$ and
            $$
            \bigvee_{\ell = 1}^{2}( \forall \gamma \in \dom(r_{\ell}))[\gamma \in \dom(p) \cap \beta \wedge r_{\ell}(\gamma) \subseteq p(\gamma)].
            $$ 
        \end{enumerate}
    \end{itemize}

    (3) If $\bfq \in \bfQ_{\bfp}$ and $\alpha \leq \lg(\bfq)$ then $\bfq \rest \alpha \in \bfQ_{\bfp}$. 

    (4) $\leq_{\bfp}$ is a partial order on $\bfQ_{\bfp}$. 

    (5) If $\bar{\bfq} = \langle \bfq_{j} \colon j < \delta \rangle$ is $\leq_{\bfp}$-increasing \underline{then} it has a $\leq_{\bfp}$-lub, $\lim(\bar{\bfq})$, 
of length $ \cup \{ \lg (\mathbf{q} _j)
  :  j < \delta \} $.  

    (6) If $\beta < \lg(\bfq)$, $\bfa = \bfa_{\bfq, \beta}$, $u \in [\cU_{\bfa, \beta}]^{\leq 2} $ 
    and $N_{u} = N_{\bfa, u}$, \underline{then}: 

    \begin{itemize}
        \item[$(\ast)$]  
        if $p \in \bbP_{\bfq}' 
       $ 
then $ q = p \rest N_{\mathbf{q}, \beta,u}$  
satisfies $ q \in N_{
u }$ and $
 q \le_{\mathbb{P} _ \mathbf{q} } p $  
\underline{where}
$ q 
  $ is defined by:  
\begin{itemize}  
\item[$ \bullet _1$]  
   $ \dom(q) = \dom(p) \cap N_{ 
u} 
\cap \beta   
$ 
\item[$ \bullet _2$] 
   If $ \gamma \in \dom(q)$ then $ q(\gamma ) = p(\gamma ) \cap N_{ u}$. 
\end{itemize}

    \end{itemize}

    (7) If (A) then (B), where: 


        \begin{enumerate}
            \item[(A)] 

            \begin{enumerate}
                \item[(a)] $i_{\ast} < \lambda$, 

                \item[(b)] $p_{i} \in \bbP_{\bfq}'$ for $i < i_{\ast}$,

                \item[(c)] if $i < j < i_{\ast}$, \underline{then} $p_{i}$ and $p_{j}$ are essentially comparable, i.e.: 

                \begin{itemize}
                    \item if $\beta \in \dom(p_{i}) \cap \dom(p_{j})$ \underline{then} $p_{i}(\beta) \subseteq p_{j}(\beta)$ or $p_{j}(\beta) \subseteq p_{i}(\beta)$. 
                \end{itemize}

                \item[(d)] $\bar{p} = \langle p_{i} \colon i < i_{\ast} \rangle$. 
            \end{enumerate}

            \item[(B)] $\bar{p}$ has a lub $p$ called $\lim(\bar{p})$ or $\lim(\{ p_{i} \colon i < i_{\ast} \})$ defined by: 

             \begin{itemize}
                    \item $\dom(p) = \bigcup \{ \dom(p_{i}) \colon i < i_{\ast} \}$, 

                    \item if $\beta \in \dom(p)$, then $$p(\beta) = \bigcup \{ p_{i}(\beta) \colon i < i_{\ast} \text{ satisfying } \beta \in \dom(p_{i}) \}.$$
                \end{itemize}
\end{enumerate}
\end{claim}

\begin{PROOF}{\ref{d11}}
    Part (2) is crucial but easy to verify. Parts (0), (1), (3), and (4) are also easy.
    
    (5) For this, define $\bfq \coloneqq \lim(\bar{\bfq})$ naturally, but we elaborate.

    \begin{enumerate}
        \item[$(\ast)$]

        \begin{enumerate}
            \item[(a)] $\lg(\bfq) = \bigcup \{ \lg(\bfq_{i}) \colon i < \delta \}$,

            \item[(b)] if $i < \delta$ and $\alpha \leq \lg(\bfq_{i})$, then $(\bbP_{\bfq, \alpha}, \bbP_{\bfq, \alpha}') = (\bbP_{\bfq_{i}, \alpha}, \bbP_{\bfq_{i}, \alpha}')$, 

            \item[(c)] if $i < \delta$ and $\beta < \lg(\bfq_{i})$, then $(\name{\bbQ}_{\bfq, \beta}, \bfa_{\bfq, \beta}, \bfb_{\bfq, \beta}) = (\name{\bbQ}_{\bfq_{i}, \beta}, \bfa_{\bfq_{i}, \beta}, \bfb_{\bfq_{i}, \beta})$, 

            \item[(d)] $(\bbP_{\bfq, \lg(\bfq)}, \bbP_{\bfq, \lg(\bfq)}')$ is $
            ( 
            \bigcup \{ \bbP_{\bfq_{i}} \colon i < \delta \}, \bigcup \{ 
            \bbP'_{\bfq_i}  
            \colon i < \delta \})$ when $\cf(\delta) \geq \lambda$, 

            \item[(e)] if $\cf(\delta) < \lambda$, then $(\bbP_{\bfq, \lg(\bfq)}, \bbP'_{\bfq, \lg(\bfq)})$ are defined as inverse limit. Then,
            \begin{itemize}
                \item $\bbP_{\bfq}' \coloneqq \bbP_{\bfq, \lg(\bfq)}'$ is dense in $\bbP_{\bfq}$ because by Definition~\ref{d5}(3), for each $\beta < \lg(\bfq_{j})$ with $j < \delta$, $\bbQ_{\bfb[\beta, \bfq_{j}]}$ is closed under increasing unions of length $< \lambda$.
            \end{itemize}
        \end{enumerate}
    \end{enumerate}

    Recalling that in Definition~\ref{d8}(1)(c), we use $\beta$ and not $\alpha$, ``$\bbP_{\bfq}$ satisfies the $\lambda^{+}$-cc'' is not required for proving \ref{d11}~(5), only ``if $\beta < \lg(\bfq)$ then $\bbP_{\bfq, \beta}$ satisfies the $\lambda^{+}$-cc'', which is clear.
Note that  even though we formally do not need it here,
the chain condition of $ \mathbb{P} _  \mathbf{q} $ 
will be proved in claim \ref{d14}.

(6)  
Note that: 
\begin{enumerate} 
\item[(a)] If $ \gamma \in 
 \dom(q) 
 $ then 
   $ \gamma \in N_ u$ and $ q(\gamma )
    \subseteq N_ u $, 
\item[(b)]  As $ \dom(q) $ and $ q(\gamma)$ 
   for $ \gamma \in \dom(q)$  has cardinality $ < \lambda $  and $ [N_u]^{< \lambda } \subseteq N_u$ 
   so recalling clause (a)  
   obviously $ q \in N_u$.

\item[(c)] To prove $ q $  is in 
$ \bbP_\mathbf{q}'$  we need, for $ \gamma \in \dom(q)$  and $ \zeta _1 < \zeta _2 $ 
from $  q( \gamma ) \subseteq  
{\mathscr U } _ \gamma $ to verify the condition in \ref{d11}(2)(c).

\item[(d)] But as 
$ \gamma \in N_u $ hence  $ \mathbf{q} \upharpoonright (\gamma +1 )$ and 
$ \zeta_1, \zeta_2$ belong to $ N_u$, also 
$ N_{\mathbf{q}, \gamma , \{ \zeta_1\} }$,
$ N_{\mathbf{q}, \gamma , \{ \zeta_2 \}  }$,
$ N_{\mathbf{q}, \gamma , \{ \zeta_1, \zeta_2 \} }$
 belong to $ N_u $  hence are included in it so we can finish easily.
\end{enumerate} 


    (7) Follows by
    our definitions. 
\end{PROOF}

We now arrive to the

\begin{cruc}\label{d14}
    If $\bfq \in \bfQ_{\bfp}$ \underline{then} $\bbP_{\bfq}$ satisfies $\lambda^{+}$-cc. Moreover $\bbP_{\bfq}$ is $\lambda^{+}$-Knaster. 
\end{cruc}

\begin{PROOF}{\ref{d14}}
    It suffices, by \ref{d8}(1)(e), to prove that $\bbP_{\bfq}' =  \bbP_{\bfq, \lg(\bfq)}'$ satisfies the $\lambda^{+}$-cc, so assume: 

    \begin{enumerate}
        \item[$(\ast)_{1}$] 

        \begin{enumerate}
            \item[(a)] Let $ \bar{p} = \langle p_{\xi} \colon \xi < \lambda^{+} \rangle$, where $p_{\xi} \in \bbP_{\bfq}'$, 

            \item[(b)] it suffices to prove that for some $\zeta < \xi < \lambda^{+}$, $p_{\zeta}$ and $p_{\xi}$ are compatible.
        \end{enumerate}
    \end{enumerate}

    [Why? By the definitions.]

    \begin{enumerate}
        \item[$(\ast)_{2}$] For some stationary set $S \subseteq \mathrm{cof}(\lambda) \cap \lambda^{+}$, we have: 

        \begin{enumerate}
            \item[$\bullet_{1}$] $\langle \dom(p_{\xi}) \colon \xi \in S \rangle$ is a $\Delta$-system with heart $w_{\ast} \in [\lg(\bfq)]^{< \lambda}$, and 

            \item[$\bullet_{2}$] if $\beta \in w_{\ast}$ then $\langle p_{\xi}(\beta) \colon \xi \in S \rangle$ is a $\Delta$-system. 

        \end{enumerate}
    \end{enumerate}

    [Why? By the Delta system lemma, the proof using Fodor's lemma recalling $\lambda = \lambda^{< \lambda}$.]

    \begin{enumerate}
        \item[$(\ast)_{3}$] Without loss of generality, $\langle p_{\xi} \colon \xi \in S \rangle$ are pairwise isomorphic (see Definition~\ref{d8}(4)).
    \end{enumerate}

    [Why? Easy  because for every $ \mathbf{a} , u $  the model $ N_{\mathbf{a}, u }$  
    has cardinality $ \lambda $.]

    \begin{enumerate}
        \item[$(\ast)_{4}$] 
For $ \gamma < \beta $ from $w_*$, we have:  
\begin{itemize}
   \item[$ \bullet _1$] Let 
   $ {\mathscr W } _{\beta}, u_{\beta, \gamma}$    
   be as in \ref{d8}(1)(d)$\bullet _5$.
   \item[$\bullet _2$]   \Wilog, $ u_{\beta, \gamma }$ 
     is disjoint to 
   $ N_{\mathbf{q} , \beta, \{ \zeta \} } \setminus N_{\mathbf{q}, \beta, \emptyset } \cap \mu $ for every 
   $ \zeta  \in {\mathscr U } _ \beta $  
   and  is  
   disjoint to 
    $ N_{\mathbf{q} , \beta, \{ \varepsilon, \zeta \} } \setminus N_{\mathbf{q}, \beta, \emptyset } \cap \mu $ for every $ \varepsilon < \zeta   $ 
      from $ {\mathscr U } _ \beta $. 
\end{itemize}
    [Why? As for any $\gamma < \beta$ from $w_{\ast}$ we have to omit from $\cU_{\beta}$ at most two elements and $w_{\ast}$ has cardinality $< \lambda$.]
\end{enumerate} 

        \begin{enumerate} 
   \item[$ (\ast)_5$]      
        We fix $\xi(1) \neq \xi(2)$ from $S$ and we shall prove that $p_{\xi(1)}$ and $p_{\xi(2)}$ have a common upper bound; this suffices for proving the Crucial Claim~\ref{d14}. 
    \end{enumerate}

\begin{enumerate}  
\item[$ (\ast)_6$]  For $ \beta \in w_*$: 
\begin{itemize}  
\item[(a)] 
       for $\ell \in \{ 1, 2\}$, consider the sequence $\langle \alpha^\beta _{\xi(\ell), \varp} \colon \varp < \varp_{\beta } \rangle$ listing  the set $p_{\xi(\ell)}(\beta)$ in increasing order 

\item[(b)]        
         Why $ \varepsilon _ \beta $ and not $ \varepsilon _{ 
        \beta, {\ell} } $? as  
       the two sequences have the same length because $p_{\xi(1)},  p_{\xi(2)}$ are isomorphic, see Definition~\ref{d8}(4)~$\bullet_{1}$.  
\item[(c)]        
       Let $\cS_\beta  \coloneqq \{ \varp < \varp_{\beta } \colon \alpha^\beta _{\xi(1), \varp}  \neq \alpha^\beta _{\xi(2), \varp}\}$,
       
\item[(d)] 
so by Definition \ref{d8}~(4)~$\bullet_{2}$ the sets $\{ \alpha^\beta _{\xi(1), \varp} \colon \varp \in \cS_{\beta} \}$, $\{ \alpha^\beta _{\xi(2), \varp} \colon \varp \in \cS_{\beta} \}$ are disjoint and disjoint to $\{ \alpha^\beta _{\xi(1), \varp} \colon \varp \in \varp_{\beta } \setminus \cS_\beta  \} = \{ \alpha^\beta _{\xi(2), \varp} \colon \varp \in \varp_{\beta } \setminus \cS _\beta \}$. 
\end{itemize}
\end{enumerate} 

    Let $\bar{\beta} = \langle \beta_{i} \colon i \leq i_{\ast} \rangle$ list the closure of $\{ \alpha, \alpha+ 1 \colon \alpha \in w_{\ast} \} \cup \{0, \lg(\bfq)\}$ in increasing order, so necessarily $i_{\ast} < \lambda$ and clearly it suffices: 

    \begin{enumerate}
        \item[$(\ast)_{7}$]  
        To choose $q_{i} \in \bbP_{\bfq, \beta_{i}}'$ a common upper bound of $\{ p_{\xi(1)} \, {\rest} \, \beta_{i},  p_{\xi(2)} \, {\rest} \, \beta_{i} \}$ increasing with $i \leq i_{\ast}$ by induction on $i \leq i_{\ast}$ 
        such that:  

\begin{enumerate}   
\item[$(\ast)$] If $\beta \in w_* \setminus \{ \beta_{j} : j < i\}$ and  $\zeta(1), \zeta(2)$ are from $ \mathscr{S}_{\beta}$
then:

\begin{enumerate} 
   \item[$\bullet _1 $] $\dom(q_j ) \cap N_{\beta, \{ \alpha _{\xi(1), \zeta(1)}, \alpha_{\xi(2), \zeta(2)}\}}$ is a subset of 
      \[
      N_{\beta, \{\alpha _{\xi (1), \zeta (1)} \} }
         \cup N_{\beta, \{\alpha _{\xi (2), \zeta (2)} \} }
         \cup N_{\beta, \emptyset},
      \] 

   \item[$ \bullet  _2$] if $ {\ell} = 1,2 $  and 
   $ \gamma \in \dom(q_j) \cap  N_{\beta, \{ \alpha  _{\xi ( {\ell} ), \zeta ( {\ell} )}\} } $ 
     then $ q_{i}( \gamma ) = p_ {\xi ({\ell} )}( \gamma )$  or
       $ \gamma \in N_{\beta, \emptyset }$
\end{enumerate} 
\end{enumerate} 



\end{enumerate}   
        \noindent
        Let us carry the induction. 

        \noindent
        \underline{Case 1}: $i = 0$. Clearly, this case is trivial,
      letting $ q_0 = \emptyset $.  

        \noindent
        \underline{Case 2}: $i$ is a limit ordinal. 
        
        In this case, let $q_{i} \coloneqq \lim \langle q_{j} \colon j < i \rangle$, so by Claim~\ref{d11}(1), $q_{i}$ is well-defined and is as required by the definition of the order 
        and satisfies $(\ast)_7$.  

        \noindent
        \underline{Case 3}: $i = j +1$ and $\beta_{j} \notin w_{\ast}$.
        
        In this case, $\dom(p_{\xi(1)}) \cap \dom(p_{\xi(2)}) \cap \beta_{i} \subseteq \beta_{j}$, hence the condition 
        $$
        q_{i} \coloneqq q_{j} \cup \left( p_{\xi(1)} \rest [\beta_{j}, \beta_{i}] 
        ) 
        \cup (p_{\xi(2)} \rest [\beta_{j}, \beta_{i}]) \right)
        $$
        is as promised. 

        \noindent
        \underline{Case 4}: $i = j +1$ and $\beta_{j} \in w_{\ast}$.
        
        By the choice of $\bar{\beta}$, clearly $\beta_{i} = \beta_{j} + 1$ 
        and let $ {\mathscr S } = {\mathscr S }_{\beta_{j}} $. 
        

        Recalling \ref{d8}(1)(d)
        and \ref{z11}(b)$(\bullet_{8})$, 
        we have: 

        \begin{enumerate}
            \item[$(\ast)_{8 
            }$]
            $\bfa_{\beta_{j}} = \bfa_{\bfq, \beta_{j}}$ determine: 

            \begin{enumerate}
                \item[(a)] $\bar{\pi}_{\beta_{j}} = \langle \pi_{u, v} \colon u, v \in [\cU_{\beta_{j}}]^{\leq 2} \text{ and }  \vert u \vert = \vert v \vert \rangle$, 

                \item[(b)] $\bar{N} _{\beta _j } 
                = \langle N_{u} \colon u \in [\cU_{\beta_{j}}]^{\leq 2} \rangle$,

                \item[(c)] for $\varp(1), \varp(2) \in \cS$, let: 

                \begin{itemize}
                    \item $v[\varp(1), \varp(2)] = \{ \alpha_{\xi(1), \varp(1)}, \alpha_{\xi(1), \varepsilon 
                    (2)}\}$, and 

                    \item $u[\varp(1), \varp(2)] = \{ \alpha_{\xi(1), \varp(1)}, \alpha_{\xi(2), \varepsilon 
                    (2)} \}$.
                \end{itemize}

                \item[(d)] for $\varp \in \cS$, let $v[\varp] = \{ \alpha_{\xi(1), \varp} \}$ and $u[\varp] = \{ \alpha_{\xi(2), \varp} \}$, 

                \item[(e)] $\bar{\iota} = \bar{\iota}_{\beta_{j}}^{\ast}$, see \ref{d8}~(1)~(d)~$\bullet_{1}  $. 

                \item[(f)] $\gamma_{j} = \xi_{\bfq}(\beta_{j})$; see \ref{d8}(1)(d) 
                $\bullet_{3} $. 
            \end{enumerate}
        \end{enumerate}

        We shall now define $p_{\varp(1), \varp(2)}$ for $\varp(1), \varp(2) \in \cS$ such that: 

        \begin{enumerate}
            \item[$(\ast)_{9 
            }$]

            \begin{enumerate}
                \item[(a)] $p_{\varp(1), \varp(2)} \in \bbP_{\gamma_{j}} \cap N_{ u[\varp(1), \varp(2)]}$, 
            hence $\dom(p_{\varp(1), \varp(2)}) \subseteq \gamma_{j} \cap N_{u[\varp(1), \varp(2)]
            }$, 

                \item[(b)] if $\varp(1) = \varp(2)$, then $p_{\varp(1), \varp(2)} \rest (\gamma_{j} \cap N_{v[\varp(1)]})$, $p_{\xi(1)} \rest N_{v[\varp(1)]}$ are essentially comparable; see \ref{d11}(7)(A)(c), moreover the first is $\leq_{\bbP_{\bfq}}$-above the second, 

                \item[(c)] if $\varp(1) = \varp(2)$, then  $p_{\varp(1), \varp(2)} \rest (\gamma_{j} \cap N_{
                u 
                [\varp(2)]})$, $p_{\xi(2)} \rest N_{
                u 
                [\varp(2)]}$ are essentially comparable, moreover the first is $\leq_{\bbP_{\bfq}}$-above the second, 

                \item[(d)] $p_{\varp(1), \varp(2)}$ satisfies \ref{d8}(1)(e)$\bullet$ with $(\gamma_{j}, \varp(1), \varp(2))$ here standing for $(\beta, \zeta_{
                1 
                }, \zeta_{2})$ there, 


                \item[(e)] $\{ q_{j} \rest N_{\emptyset} \} \cup \{ p_{\varp(1), \varp(2)} \rest N_{\emptyset} \colon \varp(1), \varp(2) \in \cS \}$ are pairwise essentially comparable, 

                \item[(f)] if $\varp(1) \neq \varp(2)$ then $p_{\varp(1), \varp(2)} \rest N_{\{  \alpha_{\varp(\ell)} \}} \leq p_{\xi(\ell)} \rest N_{\{ \alpha_{\varp(\ell)} \}}$ 
             for $  {\ell} = 1,2$. 
             
          \item[(g)]     if $ {\mathscr S } _* \subseteq {\mathscr S } \times {\mathscr S } $  then the lub $ q_{{\mathscr S } _*}  $ of 
          $ \{ q_j [N_{u[\varp(1), \varp(2)}] \colon \varepsilon(1),  \varepsilon (2) \in \mathscr{S}_{\ast} \}$ 
   satisfies the condition in $ (*)_7$. 
            \end{enumerate}
        \end{enumerate}

        We have to show two things: $\boxplus_{1}$ and $\boxplus_{2}$. The first says we can choose them (the $p_{\varp(1), \varp(2)}$-s), the second that this is enough. 

        \begin{enumerate}
            \item[$\boxplus_{1}$] we can choose $p_{\varp(1), \varp(2)}$ for $\varp(1), \varp(2) \in \cS$ as required in $(\ast)_{7}$. 
        \end{enumerate}

        We consider two possible cases: 

        \noindent
        \underline{Case 4.1}: $\varp(1) \neq \varp(2)$. 

        Let $p_{\varp(1), \varp(2)} = \pi(p_{\xi(1)} \rest N_{ v[ \varp(1), \varp(2)]})$, where $\pi = \pi_{u[\varp(1), \varp(2)], v[\varp(1), \varp(2)]}$.


Why  is $(*)_9$ preserved? Most clauses are obvious, but 
$ (*)_9(g)$ deserve elaboration, recalling that we have to satisfy $(\ast)_{7}$. 

So let $\beta \in \mathscr{W}_{\ast} \setminus \{ \beta_{\iota} \colon \iota < i \}$, hence for some $j(\ast) < i_{\ast}$, we have $\beta = \beta_{j(\ast)}$, hence we have $ \beta _{j(*)} \ge \beta _i $  
hence $ \beta _{j(*)} >  \beta_{j} $ and we have $ {\mathscr S } _ * \subseteq {\mathscr S } \times {\mathscr S }$ and deal with $q_{{\mathscr S } _* }$. 

For this, it is enough to consider the cases:

\begin{enumerate}
    \item[$\oplus_{1}$] ${\mathscr S } _* = \{  \zeta (1), \zeta (2)\}, $ 
    where $\zeta(1) = \varp(1)$ and $\zeta(2) = \varp(2)$ hence from $ {\mathscr S }$, so $ \zeta (1)\not=  \zeta (2)$,  

    \item[$\oplus_{2}$] ${\mathscr S } _* = \{  \zeta (1), \zeta (2)\} $ 
    where $ \zeta (1)\not=  \zeta (2)$ are from ${\mathscr S } $    but $ ( \zeta (1),\zeta (2)) \not= (\varepsilon(1), \varepsilon (2))$. 
\end{enumerate}

Easy to check.


        \noindent
        \underline{Case 4.2}: $\varp(1) = \varp(2)$. 
        
        In this case, we   
         pick some sequence $\langle p_{\varp, \varp} \colon \varp \in \cS \rangle$ by choosing $p_{\varp, \varp}$ by induction on $\varp \in \cS$. Now,  
        $
        p_{\varp, \varp} \in \bbP_{\beta_{j}}' \cap N_{u[\varp(1), \varp(2)]}
        $ is such that: 
        \begin{enumerate}
            \item[$(\ast)$] 

            \begin{enumerate}
                \item[(a)] $p_{\varp, \varp}$ is $\leq_{\bbP_{\bfq, \beta_{j}}'}$-above $p_{\xi(1)} \, {\rest} \, N_{v[\varp]}$ and above the restriction $p_{\xi(2)} \rest N_{u[\varp]}$, 

                \item[(b)] $\langle p_{\zeta, \zeta} \rest N_{\emptyset} \colon \zeta \in (\varp + 1) \cap \cS \rangle$ is $\leq_{\bbP_{\beta[j]}}$-increasing, and

                \item[(c)] there are $q_{1}, q_{2}, r_{1}, r_{2}$ as in Definition~\ref{d8}(2)(c) $(\bullet_{1})$-$(\bullet_{5})$ with $\bfb_{\bfq, \beta_{j}}$ standing here for $(\bfa, p, \bar{\iota})$ there such that:
                \[\bigvee_{\ell = 1}^{2}(\forall \gamma \in \dom(r_{\ell})) \left[\gamma \in \dom(p_{\varp, \varp}) \wedge r_{\ell}(\gamma) \subseteq p_{\varp, \varp}(\gamma) \right].\]
            \end{enumerate}
        \end{enumerate}

We can choose $ p_{\varepsilon, \varepsilon }$ 
by the properties of $ \mathbf{b} _{\beta _j}$


    Having defined all the $p_{\varp(1), \varp(2)}$-s we can proceed. 

    \begin{enumerate}
        \item[$\boxplus_{2}$] The following set of members of $\bbP_{\beta_{i}}$ has a common upper bound 
        $ q_*$: 

        \begin{itemize}
            \item $p_{\xi(1)}, p_{\xi(2)}$, and

            \item $p_{\varp(1), \varp(2)}$ for $\varp(1), \varp(2) \in \cS$. 
        \end{itemize}
    \end{enumerate}

    [Why? Recall Claim~\ref{d11}(2) and~\ref{d5}(1)(c)$(\bullet_{1})$ by \ref{d11}(7),
  clause (A)   there holds, in particular sub-clause 
    (A)(c). The main point is that: 

    \begin{enumerate}
        \item[$(\ast)$] $\langle N_{v[\varp(1), \varp(2)]} \cap \gamma_{j} \setminus ( N_{v[\varp(1)]} \cup N_{u[\varp(1)]} ) \colon \varp(1), \varp(2) \in \cS \rangle$ is a sequence of pairwise disjoint sets. 
    \end{enumerate}

    Why? As ``$N_{u} \cap N_{v} \subseteq N_{u \cap v}$ for $u, v \in [\cU_{\beta_{j}}]^{< 2}$  by \ref{z11}$\bullet_{7}
    $. 
    
    So $ q_* $  from $\boxplus_{2}$ 
    is a common upper bound of $p_{\xi(1)}, p_{\xi(2)}$,
    as promised. 
\end{PROOF}

\begin{remark}\label{d15}  
\begin{enumerate} 
 1)   No need so far, but we may add in $(\ast)_{4}$ of the proof of Crucial Claim~\ref{d14} the following item: 

    \begin{enumerate}
        \item[(d)] if $\beta \in w_{\ast}$ and $\langle \alpha_{\zeta, \beta, i} \colon i < \iota_{\zeta, \beta} \rangle$ lists in increasing order the members of $p_{\zeta}(\beta)$ for $\zeta \in S$, then: 

        \begin{itemize}
            \item $\langle \iota_{\zeta, \beta} \colon \zeta \in S \rangle$ is constant called $i_{\beta}$, 

            \item for $i < i_{\beta}$, the sequence $\langle \alpha_{\zeta, \beta, i} \colon \zeta \in S \rangle$ is constant or increasing, 

            \item if $i, j < i_{\beta}$ the sequence of  truth values
            $$
            \langle \text{Truth value}(\alpha_{\zeta, \beta, i} < \alpha_{\xi, \beta, j}) \colon \zeta < \xi \text{ are from } S \rangle
            $$ is constant, and

            \item if $i, j < i_{\beta}$, $\zeta \neq \xi$ are from $S$ and $\alpha_{\zeta, \beta, i} = \alpha_{\xi, \beta, j}$ then $i = j$. 
        \end{itemize}
\end{enumerate} 

2) We can make our choice of $ q_1, q_2, r_1,r_2 $ canonical, that is:

\begin{enumerate} 
\item[(A)] In \ref{d5}(2) we replace $(\bfa, p, \bar{\iota})$ by
$(\bfa, p, \bar{i}, \mathbb{F})$, where:   
\begin{itemize}
\item[$ \bullet _1$] 
    $ \mathbb{F} _{\zeta_1, \zeta_2}(q) =
    (q_1, q_2, r_1, r_2)   =
    \langle \mathbb{F} _{\zeta_1, \zeta_2, {\ell} }
(q): {\ell} = 1,2,3,4  \rangle $

\item[$ \bullet _2$] 
    if also $\zeta _3 < \zeta _4 $ 
    are   from 
    ${\mathscr U } $  then 
  $ \pi ^ \mathbf{a} _{{ \zeta _3, \zeta _4},
{ \zeta _1, \zeta _2 } } 
  \mathbb{F} _{  \zeta_1, \zeta _2,  {\ell} }
  = \mathbb{F} _{\zeta _3, \zeta _4, {\ell}},$

  where if $p\leq q \in \bbP_{\bfa} \cap N_{\bfa, \{ \varp[\bfa] \}}$ and $\zeta_{1} < \zeta_{2}$ are from $\cU$, then $\langle \mathbb{F}_{\zeta_{1}, \zeta_{1}, \ell}(p, q) \colon \ell < \mu \rangle$ is the quadruple $(q_{1}, q_{2}, r_{1}, r_{2})$ as in \ref{d5}(1)(c)($\bullet_{1}$)-($\bullet_{5}$).
  \end{itemize}

  \item[(B)]  In \ref{d5}(3) similarly 
and in \ref{d8}(1)(d)
  
  \item[(C)] In \ref{d11}(1)(d) use $ \mathbb{F} _ \beta $,

  \item[(D)]  In the proof of \ref{d14}, in
  $ (\ast)_7 \boxplus_1, $ case 4.2$(\ast)_{4.2}$ 
  we use $ \mathbb{F} 
  _{\beta _j} $,

  \item[(E)] Update the proof of \ref{d17} accordingly.
\end{enumerate} 
    \end{enumerate}
\end{remark}

\begin{claim}\label{d17}
    If (A) then (B), where: 

    \begin{enumerate}
        \item[(A)] 

        \begin{enumerate}
            \item[(a)] $\bfq \in \bfQ_{\bfp}$, 

            \item[(b)] $2 < \sigma < \lambda$,

            \item[(c)] $\name{\bfc}$ is a $\bbP_{\bfq}$-name of a function from $[\theta]^{2}$ into $\sigma$.

            \item[(d)]   $ p_* \in \mathbb{P} _  \mathbf{q} $.  
        \end{enumerate}

        \item[(B)] There is some $\bfb \in \bfA^{+}$ such that $\bbP_{\bfb} = \bbP_{\bfq}'$ and $\name{\bfc}_{\bfb} = \name{\bfc}$  
         and $ p_* \le _{\mathbb{P} _ \mathbf{q} } p_ \mathbf{b} $. 
    \end{enumerate}
\end{claim}

\begin{PROOF}{\ref{d17}}
    Recalling Hypothesis~\ref{d2}(b), on the one hand, it is clear how to choose $\bfa \in \bfA$ such that $\bbP_{\bfa} = \bbP_{\bfq}'$ and $\name{\bfc}_{\bfa} = \name{\bfc}$. On the other hand, the choice of $p_{\bfb}$ and $\bar{\iota}_{\bfb}$ is similar to the proof of \cite[2.1]{Sh:276}. We now elaborate. 

    First, we can find $\bfa$ such that: 

    \begin{enumerate}
        \item[$(\ast)_{\bfa}^{1}$] 

        \begin{enumerate}
            \item[(a)] $\bfa \in \bfA$, 

            \item[(b)] $\bbP_{\bfa} = \bbP_{\bfq}'$, 

            \item[(c)] $\gamma = \lg(\bfq)$, 

            \item[(d)] $\name{\bfc}_{\bfa} = \name{\bfc}$.
        \end{enumerate}
    \end{enumerate}
 
    Why can we find? Because we have chosen $\bbP_{\bfa}$ as in $(\ast)_{\bfa}^{1}$(b),  it is $\lambda^{+}$-cc by Claim~\ref{d14}; also $\gamma, \name{\bfc}_{\bfa}$ are as is required in Definition~\ref{d5}.  Lastly we can choose $(\cU_{\bfa}, \bar{N})$ as required because $\theta \to_{{\mathrm{sq}}} (\partial)_{
   \lambda } 
^{\lambda, 2}$ holds by Hypothesis~\ref{d2} clause (b) 
    and \ref{z11} in particular 
    clause $ (b) \bullet _{10}$. 
    
    We are left with choosing some appropriate $(p, \bar{\iota})$ and then 
 let 
    $\bfb = (\bfa, p, \bar{\iota})$. 
    Let
    \begin{multline*}
            Y \coloneqq \{ (q_{1}, q_{2}) \colon q_{1}, q_{2} \in \bbP_{\bfa}' \cap N_{\bfa, \{ \varp[\bfa] \}} \text{ are above } p_{\ast} \text{ and, }\\
            q_{1} \rest (N_{\bfa, \emptyset} \cap \lg(\bfq)) = q_{2} \rest (N_{\bfa, \emptyset} \cap \lg(\bfq)) \},
    \end{multline*}
    and let $\leq_{Y}$ be the following two place relation on $Y$:

    \begin{enumerate}
        \item[$(\ast)_{2}$] $(p_{1}, p_{2}) \leq_{Y} (q_{1}, q_{2})$ \underline{iff}: 

        \begin{enumerate}
            \item[(a)] $(p_{1}, p_{2}) \in Y$ and $(q_{1}, q_{2}) \in Y$, 

            \item[(b)] $p_{1} \leq_{\bbP_{\bfq}'} q_{1}$ and $p_{2} \leq_{\bbP_{\bfq}'} q_{2}$. 
        \end{enumerate}
    \end{enumerate}

    Clearly, 

    \begin{enumerate}
        \item[$(\ast)_{3}$] $(Y, \leq_{Y})$ is a $(< \lambda)$-complete partial order.
    \end{enumerate}

    [Why? Recalling~\ref{d11}(1).]

    \begin{enumerate}
        \item[$(\ast)_{4}$] For $(p_{1}, p_{2}) \in Y$, let
        
        \begin{enumerate}
            \item[(a)]  $\solv(p_{1}, p_{2})$ be the set of pairs $(\iota_{0}, \iota_{1})$ such that for any $\zeta_{1} < \zeta_{2}$ from $\cU_{\bfa}$, there are $r_{1}, r_{2}$ such that for $\ell = 1, 2$ clauses $\bullet_{2}$-$\bullet_{5}$ of Definition~\ref{d5}(2)(c)
            hold.  

            \item[(b)] $\solv^{+}(p_{1}, p_{2}) \coloneqq \bigcap \{ \solv(q_{1}, q_{2}) \colon (p_{1}, p_{2}) \leq_{Y} (q_{1}, q_{2}) \in Y \}$.
        \end{enumerate}
    \end{enumerate}

    \begin{enumerate}
        \item[$(\ast)_{5}$] 

        \begin{enumerate}
            \item[(a)]\label{a} if $(p_{1}, p_{2}) \leq_{Y} (q_{1}, q_{2})$ then: 
            $$
            \solv(p_{1}, p_{2}) \supseteq \solv(q_{1}, q_{2}) \supseteq \solv^{+}(q_{1}, q_{2}) \supseteq \solv^{+}(p_{1}, p_{2}),
            $$ 

            \item[(b)] if $(p_{1}, p_{2}) \in Y$ then $\solv(p_{1}, p_{2}) \neq \emptyset$. 
        \end{enumerate}
    \end{enumerate}

    [Why? The first inclusion in Clause (a) holds because $\leq_{\bbP_{\bfq}}$ is transitive. The other inclusions are clear, and Clause (b) is easy too.]

    \begin{enumerate}
        \item[$(\ast)_{6}$] If $(p_{1}, p_{2}) \in Y$ then for some $(q_{1}, q_{2})$ and $\bar{\iota}$, we have: 

        \begin{enumerate}
            \item[(a)] $(p_{1}, p_{2}) \leq_{Y} (q_{1}, q_{2}) \in Y$, 

            \item[(b)] if $(q_{1}, q_{2}) \leq_{Y} (q_{1}', q_{2}')$ \underline{then} $\bar{\iota} \in \solv(q_{1}', q_{2}')$, moreover, $\solv(q_{1}, q_{2}) =  \solv(q_{2}', q_{2}') = \solv^{+}(q_{1}', q_{2}') = \solv^{+}(q_{1}, q_{2})$. 
        \end{enumerate}
    \end{enumerate}

    [Why? Recalling $\sigma < \lambda$, hence $\vert \sigma \times \sigma \vert < \lambda$ and $(Y, \leq_{Y})$ is $\lambda$-complete by $(\ast)_{3}$.]

    \begin{enumerate}
        \item[$(\ast)_{7}$] For $p \in \bbP_{\bfa}' \cap N_{\bfa, \{ \varp[\bfa] \}}$, let $\solv(p)$ be the set of $\bar{\iota} \in \sigma \times \sigma$ such that there is $(q_{1}, q_{2})$ such that: 

        \begin{enumerate}
            \item[$\bullet_{1}$] $p \leq_{\bbP_{\bfq}} q_{1}$, $p \leq_{\bbP_{\bfq}} q_{2}$ and

            \item[$\bullet_{2}$] $(q_{1}, q_{2}) \in Y$, 

            \item[$\bullet_{3}$] $\bar{\iota} \in \solv^{+}(q_{1}, q_{2})$, 

            \item[$\bullet_{4}$] $\solv(q_{1}, q_{2}) = \solv^{+}(q_{1}, q_{2})$. 
        \end{enumerate}
    \end{enumerate}

    \begin{enumerate}
        \item[$(\ast)_{8}$] 

        \begin{enumerate}
            \item[(a)] if $p \in \bbP_{\bfa}' \cap N_{\bfa, \{ \varp[\bfa] \}}$ then $\solv(p) \neq \emptyset$,

            \item[(b)] if $p \leq_{\bbP_{\bfa}'} q$ are from $\bbP_{\bfa}' \cap N_{\bfa, \{ \varp[\bfa] \}}$ then $\solv(p) \supseteq \solv(q)$, 

            \item[(c)] if $p \in \bbP_{\bfa}' \cap N_{\bfa, \{ \varp[\bfa] \}}$ then for some $q$ and $\bar{\iota}$, for every $q'$, we have 
            $ q \leq_{\bbP_{\bfq}'}  q' \wedge  q' \in \bbP_{\bfa}'
            \cap N_{\bfa, \{ \varp[\bfa] \}} \Rightarrow \bar{\iota} \in \solv(q')$. 
        \end{enumerate}
    \end{enumerate}

    [Why? Clause (a) follows by $(\ast)_{6}$, Clause (b) by the definitions, and Clause (c) holds as $\bbP_{\bfa}'$ and even $\bbP_{\bfa}' \cap N_{\bfa, \{ \varp[\bfa] \}}$ is $\lambda$-complete and $\vert \sigma \times \sigma \vert < \lambda$.]  

    Now, 
    applying  
    $(\ast)_{8}$(c) 
    to $ p_* $  
    finish the proof of~\ref{d17}.
\end{PROOF}

\begin{claim}\label{d20}
    If (A) then (B), where: 

    \begin{enumerate}
        \item[(A)] 

        \begin{enumerate}
            \item[(a)] $\bfq \in \bbQ_{\bfp}$ and $\bfq_{0} <_{\bfp} \bfq$, 

            \item[(b)] $\gamma(\bfq) < \mu$, so $\lg(\bfq) < \mu$,

            \item[(c)] $\bfb \in \bfA_{\bfp}$ and $\bbP_{\bfb} = \bbP_{\bfq_{0}}$.  
        \end{enumerate}

        \item[(B)] There exists some $\bfq_{1}$ such that: 

        \begin{enumerate}
            \item[(a)] $\bfq \leq_{\bfp} \bfq_{1}$, 

            \item[(b)] $\lg(\bfq_{1}) = \lg(\bfq) + 1$, 

            \item[(c)] $\bfb_{\lg(\bfq)}[\bfq_{1}] = \bfb$.
        \end{enumerate}
    \end{enumerate}
\end{claim}

\begin{PROOF}{\ref{d20}}
   Easy.
\end{PROOF}

Lastly, before arriving at the main conclusion, we have to prove the following. 

\begin{claim}\label{d21}\ 

    (1) Assume $\bfq \in \bfQ_{\bfp}$, $\alpha < \lg(\bfq)$ and $\bfb = \bfb_{\bfq, \alpha} = (\bfa_{\alpha}, p_{\alpha}, \bar{\iota}_{\alpha}) = (\bfa, p, \bar{\iota})$, then: 

    \begin{itemize}
        \item $\Vdash_{\bbP_{\bfq, \alpha + 1}}$``$
        \mathcal{V}
        _{\name{\bbQ}_{\bfb}} \in [\cU_{\bfa_{\alpha}}]^{\partial}$ and for every $\alpha \neq \beta \in \mathcal{V}_{\name{\bbQ}_{\bfb}}$,  $\name{\bfc}_{\bfa_{\alpha}}\{ \alpha, \beta \} \in \{ \iota_{1}, \iota_{2} \}$''. 
    \end{itemize}

    (2) If $\bfb = (\bfa, p, \bar{\iota}) \in \bfA^{+}$, $\cf(\partial) > \lambda$, and in 
    $  \bfV   
    ^{\bbP_{\bfa}}$, $\bbQ_{\bfb}$ satisfies the $\lambda^{+}$-cc, \underline{then} for some $p \in \bbQ_{\bfb} \cap \bbP_{\bfa} \cap N_{\bfa, \{ \varp[\bfa] \}}$ we have\footnote{We may omit $p$ but it does not matter.} 
    $p \Vdash_{\name{\bbQ}_{\bfb}}$``$\mathcal{V}_{\name{\bbQ}_{\bfb}} \in [\cU_{\bfa}]^{\partial}$ and for every $\alpha \neq \beta \in \mathcal{V}_{\name{\bbQ}_{\bfb}}$, 
    $\bfc_{\bfa  
    }
    \{ \alpha, \beta \} \in \{ \iota_{1}, \iota_{2} \}$''.  
\end{claim}

\begin{PROOF}{\ref{d21}}
    (1) The second phrase in both conclusion holds by the definitions of $\name{\bbQ}_{\bfb}$. 

    By the proof of ``$\bbP_{\bfq}$ satisfies the $\lambda^{+}$-cc'', we can show for $\varp < \partial$, the density of the set 
    $$
    \cI_{\varp} \coloneqq \{ p \in \bbP_{\bfq}' \colon \alpha \in \dom(p) \text{ and there is } \beta \in p(\alpha) \text{ such that } \varp < \otp(\cU_{\bfa_{\alpha}} \cap \beta)   \}.
    $$

    (2) Easily, for every $\beta \in \cU_{\bfa}$ we can choose $p_{\beta}^{0} = \{ \beta \}$, $q_{\beta} = \{ (p, p_{\beta}^{0}) \}$. Clearly, $q_{\beta} \in \bbP_{\alpha} \ast \name{\bbQ}_{\bfb}$ for $\beta \in \cU_{\bfa}$. So by the $\lambda^{+}$-cc for some $\beta \in \cU_{\bfa}$, $q_{\beta} \Vdash$``$\{ \varp \in \cU_{\bfa} \colon q_{\varp} \in \name{\bbQ}_{\bfb} \} \in [\cU_{\bfa}]^{\partial}$; well assuming $\cf(\partial) > \lambda$. 
\end{PROOF}

\begin{conclusion}\label{d23}
    There exists a forcing notion $\bbP$ satisfying the following conditions: 

    \begin{enumerate}
        \item[(a)] $\bbP$ is $\lambda^{+}$-cc of cardinality $\mu$. 
        
        \item[(b)] $\bbP$ is $({<} \, \lambda)$-complete; hence, it collapses no cardinals, changes no cofinalities, and preserves cardinal arithmetic outside the interval $[\lambda, \mu)$.

        \item[(c)]  
    $\Vdash_{\bbP}$``$2^\lambda = \mu "$.  
    
        \item[(d)]  
        $\Vdash_{\bbP}$``$\theta \to [\partial]_{\sigma, 2}^{2}$'' for every $\sigma \in (2, \lambda)$. 
    \end{enumerate}
\end{conclusion}

\begin{PROOF}{\ref{d23}}
    Choose a $\leq_{\bfp}$-increasing continuous sequence $\langle \bfq_{\alpha} \colon \alpha < \mu \rangle \in {}^{\mu}({\bfQ_{\bfp}})$ such that $\lg(\bfq_{\alpha}) = \alpha$, $\bbP_{\bfq_{\alpha}}$ has cardinality $\leq (\vert \alpha \vert + \lambda)^{< \lambda}$ and,

    \begin{itemize}
        \item if $\alpha < \mu$ and $\Vdash_{\bbP_{\bfq_{\alpha}}}$``$\name{\bfc} \colon [\theta]^{2} \to \sigma$'', then for unboundedly many $\beta \in [\alpha, \mu)$, $\name{\bfc}_{\bfq_{\beta + 1, \beta}} = \name{\bfc}$.
    \end{itemize}
    
    The existence of $\bfb_{\beta}[\bfq_{\beta + 1}]$ with $\name{\bfc}[\bfb_{\beta}[\bfq_{\beta + 1}]] = \name{\bfc}$ as required holds by Claim~\ref{d17} and Claim~\ref{d20} below. 

    Clearly $\bigcup \{ \bbP_{\bfq_{\beta}} \colon \beta < \mu \}$ is a forcing notion as is required. 
\end{PROOF}

Conclusion~\ref{d23} is meaningful because:

\begin{fact}\label{d26}
    Assume that $\lambda = \lambda^{ < \lambda} < \partial < \theta < \mu = \mu^{\theta}$, and $[\alpha < \mu \Rightarrow \vert \alpha \vert^{\lambda} < \mu]$, $\theta > \beth_{2}(\kappa)$ and $\partial = \kappa^{+}$, $\kappa = \kappa^{\lambda}$. \underline{Then} the demands in Hypothesis~\ref{d2} hold. 
\end{fact}

\begin{remark}\label{d29}
    To justify the assumption, notice that: 

    \begin{enumerate}
        \item[(A)] Omitting $\kappa = \kappa^{\lambda}$ does not help. 

        \item[(B)] $\theta \to_{\mathrm{sq}} (\partial)_{\partial}^{2 \leq \lambda}$ implies $\theta \to (\partial)_{2^{\partial}}^{2}$, hence necessarily $\theta > 2^{2^{\partial}}$. 
    \end{enumerate}

    With stronger lower bound on $\theta$, see~\cite{Sh:289}. 

    The main point is proving $\theta \to_{\mathrm{sq}} (\partial)_{\partial 
    }^{\leq \lambda, 2}$. For this, see~\cite{Sh:289}, $\theta = \beth_{m}(\partial)$ for some small $m$ suffice, on this the bounds in \ref{d23} depends; we intend to return to this in  \cite{Sh:F2407}. Anyhow just $\theta < \partial^{+ \omega}$ and $\GCH$ in $[\partial, \partial^{+ \omega}]$ would suffice for me.   
\end{remark}

\begin{PROOF}{\ref{d26}}
    The point is to prove $\theta \to_{\mathrm{sq}} 
    (\partial 
    )_{\partial
    }^{\lambda, 2}$. Let $\cB$ be as in \ref{z11}(a), $\partial_{1} = 2^{\kappa}$, $\partial_{2} = 2^{\partial_{1}}$, and $\theta > \partial_{2}$. 

    Let $\chi > 2^{\mu}$, and  $\gC_{\ast}$ be an expansion of  $(\cH(\chi), \in, <_{\chi}^{\ast}, \cB)$ with vocabulary of cardinality $\lambda$ such that for any finite set $u \subseteq \cH(\chi)$, the Skolem hull of $u$, $N_{u}  \coloneqq \mathrm{Sk}(u, \gC_{\ast})$ is of cardinality $\lambda$ and $\vert N_{u} \vert^{< \lambda} \subseteq N_{u}$.

    Let $\gC_{2} \prec_{\bbL_{\partial(1)^{+}, \partial(1)^+{}}} \gC_{\ast}$ be of cardinality $\partial_{2}$ such that $\partial_{2} + 1 \subseteq \gC_{2}$. Let $\beta_{1} \coloneqq \min(\theta \setminus \gC_{2})$. Similarly, choose $\gC_{1} \prec_{\bbL_{\partial, \partial}} \gC_{\ast}$ of cardinality $\partial_{1}$ such that $\partial_{1} + 1 \subseteq \gC_{1}$ and $\{ \gC_{2}, \beta_{0} \} \subseteq \gC_{1}$.  

    Let $\gC_{0} = \gC_{1} \cap \gC_{2}$ and choose $\beta_{0} \in \beta_{1} \cap \gC_{2} \subseteq \theta \cap \gC_{2}$ realizing the $\bbL_{\partial, \partial}$-type which $\beta_{1}$ realizes over $\gC_{0}$. 

    Now,
    
    \begin{enumerate}
        \item[$(\ast)_{1}$] choose $\alpha_{\varp} \in \gC_{0} \cap \theta$ by induction on $\varp < \partial$, such that: 

        \begin{itemize}
            \item $\alpha_{\varp}$, $\beta_{1}$ realize the same first-order type in $\gC_{\ast}$ over the set $\{ \beta_{2} \} \cup (A_{\varp} \cap \gC_{0})$, where: 
            \[A_{\varp} =  \mathrm{Sk}_{\gC}( \{ \alpha_{\zeta} \colon \zeta < \varp \} \cup \{ \beta_{1}, \beta_{0} \}).
            \]
        \end{itemize}
    \end{enumerate}
    
    \begin{enumerate}
        \item[$(\ast)_{2}$] Let $N_{\emptyset}^{\bullet} = N_{\{ \beta_{0}, \beta_{1} \}} \cap \gC_{0}$. 
    \end{enumerate}

    Note, 

    \begin{enumerate}
        \item[$(\ast)_{3}$] for $\varp < \zeta < \partial$, the following pairs realize the same type over $N_{0}^{\ast}$ in $\gC_{\ast}$: 

        \begin{itemize}
            \item[$\bullet_{1}$] $(\alpha_{\varp}, \alpha_{\zeta})$,

            \item[$\bullet_{2}$] $(\alpha_{\varp}, \beta_{0})$, 

            \item[$\bullet_{3}$] $(\alpha_{\varp}, \beta_{1})$, 

            \item[$\bullet_{4}$] $(\beta_{0}, \beta_{1})$.
        \end{itemize}
    \end{enumerate}

    [Why? For the equality of $\bullet_{1}$ and $\bullet_{2}$ note the choice of $\alpha_{\varp}$. 
    
    For the equality of $\bullet_{2}$ and $\bullet_{3}$, note the choice of $\beta_{0}$.

    For the equality of $\bullet_{3}$ and $\bullet_{4}$ note the choice of $\alpha_{\varp}$.]

    \begin{enumerate}
        \item[$(\ast)_{4}$] 

        \begin{itemize}
            \item[$\bullet_{1}$] $N_{\{\varp, \zeta\}}^{\ast} = N_{\{ \alpha_{\varp}, \alpha_{\zeta} \}}$, so for $\varp < \zeta < \partial$, 

            \item[$\bullet_{2}$] $N_{\{ \varp, \zeta \}}^{\ast} \prec \gC_{0}$, 

            \item[$\bullet_{3}$] $N_{\emptyset}^{\bullet} \prec N_{\{ \varp, \zeta \}}^{\ast}$. 
        \end{itemize}
    \end{enumerate}

    \begin{enumerate}
        \item[$(\ast)_{5}$] for $\varp < \zeta < \partial$, let $f_{\{ \varp, \zeta \}}$ be the isomorphism from $N_{\{ \varp, \zeta \}}^{\ast}$ onto $N_{\{ \beta_{0}, \beta_{1} \}}$.
    \end{enumerate}

    [Why does it exist?  by $(\ast)_{3}$.]

    \begin{enumerate}
        \item[$(\ast)_{6}$] $f_{\{ \varp, \zeta \}}$ is the identity on $N_{\emptyset}^{\bullet}$ (and $N_{\emptyset}^{\bullet} \prec N_{\{ \varp, \zeta \}}$). 
    \end{enumerate}

    [Why? By $(\ast)_{2}$.]

     \begin{enumerate}
        \item[$(\ast)_{7}$] if $\varp(0) < \zeta(0) < \partial$, $\varp(1) < \zeta(1) < \partial$ and $\{ \varp(0), \zeta(0) \} \cap \{ \varp(1), \zeta(1) \} = \emptyset,$ then 
        \[
            N_{ \{ \varp(0), \zeta(0)\}}^{\ast} \cap N_{\{\varp(1), \zeta(1) \}}^{\ast} = N_{\{ \beta_{0}, \beta_{1} \}} \cap \gC_{0} = N_{\emptyset}^{\bullet}. 
        \]
    \end{enumerate}

    [Why? The second equality holds by $(\ast)_{2}$; without loss of generality $\zeta(0) < \zeta(1)$. Now, 
    
    \begin{enumerate}
        \item[$\bullet_{1}$] $N_{\{ \varp(0), \zeta(0) \}}^{\ast} \cap N_{\{ \varp(1), \zeta(1) \}}^{\ast} = N_{\{ \varp(0), \zeta(0) \}}^{\ast} \cap N_{\{ \alpha_{\varp(1), \beta_{1}} \}}$ by the choice of $\zeta(1)$. 

        \item[$\bullet_{2}$] if $\zeta(0) < \varp(1)$ then 
        
        $N_{\{ \varp(0), \zeta(0) \}}^{\ast} \cap N_{\{ \alpha_{\varp(1), \beta_{1}} \}} = N_{\{ \varp(0), \zeta(0) \} \}}^{\ast} \cap N_{\{ \beta_{0}, \beta_{1} \}} = N_{\{ \varp(0), \zeta(0) \}} \cap N_{\emptyset}^{\bullet} = N_{\emptyset}^{\bullet}$ because the first equality follows by the choice of $\alpha_{\varp(1)}$ second equality by $(\ast)_{4}\bullet_{2}$ and $(\ast)_{2}$; the third equality by $(\ast)_{3}$. 
    \end{enumerate}

    \begin{enumerate}
        \item[$\bullet_{3}$] if $\varp(0) < \varp(1) < \zeta(0)$, then: 
        \begin{equation*}
            \begin{split}
                N_{\{ \varp(0), \zeta(0) \}}^{\ast} \cap N_{\{ \alpha_{\varp(1)}, \beta_{1} \}} &= N_{\{ \alpha_{\varp(0)}, \beta_{0} \}} \cap N_{ \{ \alpha_{\varp(1)}, \beta_{1} \}}\\
                &=
                \left(  N_{\{ \alpha_{\varp(0)}, \beta_{0} \}}^{\ast} \cap \gC_{0} \right) \cap  \left( N_{\{ \alpha_{\varp(0), \beta_{1}} \}}  \cap \gC_{0} \right)\\
                &= \left(N_{\{ \alpha_{\zeta(0), \beta_{1}} \}} \cap \gC_{0} \right) \cap \left( N_{\{ \alpha_{\zeta(0), \beta_{1}}\}} \cap \gC_{0} \right)\\
                & = \left( N_{\{ \alpha_{\zeta(0)}, \beta_{1} \}} \cap \gC_{0} \right) \cap \left(  N_{\{ \beta_{0}, \beta_{1} \}} \cap \gC_{0} \right)\\
                & = \left( N_{\{ \alpha_{\varp(0)}, \beta_{1} \}} \cap \gC_{0} \right) \cap N_{\emptyset}^{\bullet} = N_{\emptyset}^{\bullet}.
            \end{split}
        \end{equation*}
    \end{enumerate}

    [Why? The first equality holds by the choice of $\beta_{0}$. The second equality as $N_{\{ \varp(0), \zeta(0) \}} \subseteq \gC_{0}$ and the first equality. The third equality holds by the choice of $\beta_{0}.$ The fourth equality  holds by the choice of $\alpha_{\zeta(0)}$. The fifth equality holds by the choice of $N_{\emptyset}^{\bullet}$ i.e., $(\ast)_{2}$. Finally, the sixth equality holds as $N_{\{ \alpha_{\varp(1)}, \beta_{1} \}} \supseteq N_{\emptyset}^{\bullet}$ by the choice of $\alpha_{\varp(2)}$.]

    \begin{enumerate}
        \item[$\bullet_{4}$] If $\varp(1) < \varp(0)$, then:
        \begin{equation*}
            \begin{split}
                N_{\{ \varp(0), \zeta(0) \}}^{\ast} \cap N_{\{ \alpha_{\varp(1)}, \beta_{1} \}} &=
                N_{\{ \alpha_{\varp(0)}, \beta_{0} \}} \cap N_{\{ \alpha_{\varp(1)}, \beta_{1} \}}\\
                & = \left(  N_{\{ \alpha_{\varp(0)}, \beta_{0} \}}^{\ast} \cap \gC_{0}  \right) \cap \left( N_{\{ \alpha_{\varp(1), \beta_{1}} \}}^{\ast} \cap \gC_{0}  \right)\\
                &= \left( N_{\{ \alpha_{\varp(0)}, \beta_{1} \}} \cap \gC_{0} \right) \cap \left(  N_{ \{ \alpha_{\varp(1)}, \beta_{1} \} } \cap \gC_{0} \right)\\
                &= \left( N_{ \{ \beta_{0}, \beta_{1} \}}  \cap \gC \right) \cap \left( N_{ \{ \alpha_{\varp(1)}, \beta_{1} \} } \cap \gC  \right)\\
                &= N_{\emptyset}^{\ast} \cap \left(  N_{ \{ \alpha_{\varp(1)}, \beta_{1} \} } \cap   \gC_{0} \right) = N_{\emptyset}^{\bullet}.\\
            \end{split}
        \end{equation*}    
    \end{enumerate}

    [Why? The first equality holds by the choice of $\beta_{0}$. The second one holds as $N_{ \{ \varp(0), \zeta(0) \} } \subseteq \gC_{0}$ and the first equality. The third equality holds by the choice of $\beta_{0}$. The fourth equality holds by the choice of $\alpha_{\varp(0)}$. The fifth equality holds by the choice of $N_{\emptyset}^{\bullet}$, i.e., by $(\ast)_{2}$. Finally, the sixth equality holds as $N_{\{ \alpha_{\varp(1)}, \beta_{1} \}} \supseteq N_{\emptyset}^{\bullet}$ and by the choice of $\alpha_{\varp(0)}$.

    Recalling $\bullet_{1}$ and the division to cases in $\bullet_{2}, \bullet_{3}$ and $\bullet_{4}$, we are done proving $(\ast)_{6}$.] 

    \begin{enumerate}
        \item[$(\ast)_{7}$] if $\varp < \zeta(1) < \zeta(2) < \partial$, \underline{then} $N_{\{ \varp, \zeta(1) \}} \cap N_{\varp, \zeta(2)} = N_{\{ \varp \}}^{\uparrow} \coloneqq N_{\{ \alpha_{\varp}, \beta_{1} \}} \rest \gC_{0}$.
    \end{enumerate}

    [Why? By  the choice of $\alpha_{\zeta(2)}$ and  $\alpha_{\zeta(1)}$.]

    \begin{enumerate}
        \item[$(\ast)_{8}$] if $\zeta_{1} < \zeta_{2} < \varp < \partial$, then $N_{\{ \zeta_{1}, \varp \}}^{\ast} \cap N_{\{ \zeta_{2}, \varp \}}^{\ast} = N_{\{ \varp \}}^{\downarrow}$, where $N_{\{ \varp \}}^{\downarrow} \coloneqq  f_{\varp, \varp + 1}^{-1}(N_{ \{\beta_{0}, \beta_{1} \}})$. 
    \end{enumerate}

    [Why? For $\zeta < \partial$, $N_{ \{ \alpha_{\zeta}, \beta_{0} \} } \cap \gC_{0} =  N_{ \{ \alpha_{\zeta_{3}}, \beta_{1} \} } \cap \gC_{0}$ by the choice of $\beta_{0}$, and $\alpha_{\varp}, \beta_{0}$ realize the same type of $\gC_{\ast}$ over $\{ \beta_{1} \} \cup (A_{\varp} \cap \gC_{0})$.]

  \begin{enumerate}
        \item[$(\ast)_{9}$] 

        \begin{itemize}
            \item Let $N_{\{ \varp \}}^{\ast}$ be the $\mathrm{Sk}(N_{\{ \varp \}}^{\uparrow} \cup N_{\{ \varp \}}^{\downarrow}, \gC_{\ast})$, and

            \item  let $M_{\varp}^{\ast} = \left(\mathrm{Sk}(\bigcup_{\ell < 5} N_{ \{ 5\varp + \ell\}}^{\ast} \cup \{ N_{\{5\varp + m, \delta \varp + n \}} \colon m < n < 5 \}^{\ast}) \} \right)_{\ell < 5}$, 

            \item let $M_{\varp}^{+}$ be $M_{\varp}^{\ast}$ expanded by: 

            \item $c_{\ell}^{M_{\varp}^{+}} = \alpha_{5 \varp + \ell}$ for $\ell < 5$, 

            \item $p_{\ell}^{M_{\varp}^{+}} = \vert N_{ \{ 5 \varp + \ell \} }^{\ast} \vert$ for $\ell < 5$, 

            \item $P_{m, n}^{M_{\varp}^{+}} = \vert N_{ \{ 5 \varp + m, 5 \varp + n \colon m < n < 5 \}}^{\ast} \vert$.
        \end{itemize}
    \end{enumerate}

     \begin{enumerate}
        \item[$(\ast)_{10}$] There is some $\cU_{1} \in [\partial]^{\partial}$ such that: 

        \begin{itemize}
            \item $\langle M_{\varp}^{\ast} \colon \varp \in \cU \rangle$ is a $\Delta$-system with heart $N_{\emptyset}^{\ast}$, 

            \item the $M_{\varp}$ are pairwise isomorphic. 
        \end{itemize}
    \end{enumerate}

    [Why? Because $\partial = \partial_{0}$ and $\partial_{0} = (\partial_{0})^{+}$ by the $\Delta$-system lemma.]

    \begin{enumerate}
        \item[$(\ast)_{11}$] $\langle N_{u}^{\ast} \colon u \in \cU_{2} \rangle$ is a required when $\cU_{2} =  \{ 5 \varp + 2 \colon \varp \in \cU_{1} \}$ and $N_{\{ 5 \varp + 2 \}}^{\ast} = M_{\varp}^{\ast}$.  
    \end{enumerate}

    Pedantically,  $\cU_{3} = \{ \alpha_{\zeta} \colon \zeta \in \cU_{2} \}$ and $N_{ \{ \alpha_{\zeta} \colon \zeta \in u \} }^{\ast} = N_{u}^{\ast}$ for $u \in [\cU_{3}]^{\leq 2}$.
\end{PROOF}

\bibliographystyle{amsalpha}
\bibliography{shlhetal}
\end{document}